\documentclass[12pt]{article}

\usepackage{amsmath}
\usepackage{amsthm}
\usepackage{amssymb}

\theoremstyle{definition}
\newtheorem{theorem}{Theorem}[section]
\newtheorem{lemma}[theorem]{Lemma}
\newtheorem{corollary}[theorem]{Corollary}
\newtheorem{question}[theorem]{Question}
\newtheorem{definition}[theorem]{Definition}

\newtheorem{example}[theorem]{Example}
\newtheorem{remark}[theorem]{Remark}

\newcommand{\setN}{\mathbb{N}}
\newcommand{\setQ}{\mathbb{Q}}

\newcommand{\UF}{\operatorname{UF}}
\newcommand{\MF}{\operatorname{MF}}

\usepackage{xspace}
\newcommand{\define}[1]{\textit{#1}\xspace}

\newcommand{\I}{\textsc{I}\xspace}
\newcommand{\II}{\textsc{II}\xspace}

\begin{document}

\title{\textbf{Topological aspects of poset spaces}}

\author{Carl Mummert%
  \thanks{ \texttt{mummertc@marshall.edu}.  C.~Mummert was
    partially supported by a VIGRE graduate traineeship
    under NSF Grant DMS-9810759 at the Pennsylvania State University.}\\[6pt]
  and Frank Stephan\footnote{\tt fstephan@comp.nus.edu.sg.\
    \rm F.~Stephan is supported in part by NUS grant number
    R252-000-212-112.}}

\date{April 3, 2008\footnote{Accepted for publication in the
    \textit{Michigan Mathematical Journal}. An earlier version of this
    paper was released as National University of Singapore School Of
    Computing technical report TRC6/06.}\\[4pt] Revised November 13,
  2009}
\maketitle

\begin{abstract}
  \noindent
  We study two classes of spaces whose points are filters on
  partially ordered sets.  Points in MF spaces are maximal
  filters, while points in UF spaces are unbounded filters.
  We give a thorough account of the topological properties
  of these spaces.  We obtain a complete characterization of
  the class of countably based MF spaces: they are precisely
  the second-countable $T_1$ spaces with the strong Choquet
  property.  We apply this characterization to domain theory
  to characterize the class of second-countable spaces with
  a domain representation.
\end{abstract}

\section{Introduction}\label{s1}

\noindent
Recent work in mathematical
logic~\cite{Mummert-thesis,Mummert-MFspaces,MumSim} has led
to an interest in certain topological spaces formed from
filters on partially ordered sets.  This paper describes the
general topology of these poset spaces.

The results of the paper are divided as follows. In
Section~\ref{sec2}, we define two classes of spaces, MF spaces and UF
spaces. Together these spaces form the class of poset spaces. We show
that many familiar spaces are homeomorphic to poset spaces.  In
Section~\ref{sec3}, we characterize the separation properties of poset
spaces and show that any second-countable poset space is homeomorphic
to a space of the same kind formed from a countable poset. In
Section~\ref{sec4}, we show that the class of MF spaces are closed
under arbitrary topological products and that any $G_\delta$ subspace
of an MF space is again an MF space.  We show that UF spaces are
closed under taking $G_\delta$ subspaces but not closed under binary
products.  In Section~\ref{sec5}, we establish that poset spaces are
of the second Baire category and possess the strong Choquet property.
We give a characterization of the class of countably based MF spaces
as the class of second-countable $T_1$ spaces with the strong Choquet
property.  In Section~\ref{sec5a}, we apply the results of
Section~\ref{sec5} to domain theory, giving a complete
characterization of the second-countable topological spaces that have
a domain representation.  Section~\ref{sec7} contains results on the
relationship between MF spaces (not necessarily countably based) and
semi-topogenous orders.  We use semi-topogenous orders to establish a
sufficient condition for an arbitrary space to be homeomorphic to an
MF space.  In Section~\ref{sec6}, we show that every second-countable
poset space is either countable or contains a perfect closed set.

\section{Poset spaces}\label{sec2}

\noindent Our goal in this section is to define the class of
poset spaces and show that this class includes all complete
metric spaces and all locally compact Hausdorff spaces.  We
first review some basic definitions about partially ordered
sets.

A \define{poset} is a set $P$ with an reflexive,
antisymmetric, transitive relation $\preceq$. That is, the
following conditions hold for all $p$, $q$ and $r$ in $P$.
\begin{enumerate}
\item $p \preceq p$.
\item If $p \preceq q$ and $q \preceq p$ then $q =p$.
\item If $p \preceq q$ and $q \preceq r$ then $p \preceq r$.
\end{enumerate}
We write $p \prec q$ if $p \preceq q$ and $p \not = q$.  If
there is no $r$ such that $r \preceq p$ and $r \preceq q$
then we write $p \perp q$.

A \define{filter} is a subset $F$ of a poset $P$ satisfying
the following two conditions.
\begin{enumerate}
\item For every $p,q \in F$ there is an $r \in F$ such that
  $r \preceq p$ and $r \preceq q$.
\item For every $p \in F$ and $q \in P$ if $p \preceq q$ then
  $q \in F$.
\end{enumerate}
A filter $F$ is \define{unbounded} if there is no $r \in P$
such that $r \prec q$ for every $q \in F$. Furthermore, $F$ is
\define{maximal} if there is no strictly larger filter
containing $F$.  Every maximal filter is unbounded, but in
general not every unbounded filter is maximal.

For any poset $P$, we let $\UF(P)$ denote the set of
unbounded filters on~$P$ and let $\MF(P)$ denote the set of
maximal filters on $P$.  We topologize $\UF(P)$ with the
basis $\{N_p \mid p \in P\}$, where
\[
N_p = \{ F \in \UF(P) \mid p \in F\}.
\]
We give $\MF(P)$ the topology it inherits as a subset of
$\UF(P)$; when we work with spaces of maximal filters we may
write $N_p$ to denote the set of maximal filters containing
$p$.  To facilitate the exposition, we sometimes identify $p
\in P$ with the open set $N_p$ and identify a subset $U$ of
$P$ with the open set $\bigcup_{p \in U} N_p$.

A \define{UF space} is a space of the form $\UF(P)$ and an
\define{MF space} is a space of the form $\MF(P)$.  UF
spaces and MF spaces are collectively referred to as
\define{poset spaces}.  A poset space is \define{countably
  based} if it is formed from a countable poset.  It is
possible that $P$ is uncountable but $\MF(P)$ or $\UF(P)$ is
a second-countable space (an example is provided after
Theorem~\ref{s2t2}).  We will show below that every
second-countable poset space is homeomorphic to a countably
based poset space.  This result justifies our terminology.

\begin{remark}
  It is sometimes convenient to work with strict partial
  orders instead of the non-strict partial orders defined
  above. A strict partial order is a set $P$ with an
  irreflexive, transitive relation $\prec$. Every strict
  partial order $\langle P, \prec\rangle$ is canonically
  associated to non-strict partial order $\langle P,
  \preceq\rangle$ in which $p \preceq q$ if and only if $p
  \prec q$ or $p = q$, and every non-strict partial order
  arises in this way.
% Conversely, every non-strict partial
%  relation $\preceq$  order canonically determines a strict
%  partial order relation $\prec$ in which $p \prec q$ if $p
%  \preceq q$ and $p \not = q$.
  A filter on a strict partial order $\langle P,
  \prec\rangle$ is a set $F \subseteq P$ that is upward
  closed and such that if $p,q \in F$ then there is an $r
  \in F$ with $r \preceq p$ and $r \preceq q$. 
  
  It follows immediately from these definitions that if
  $\langle P, \prec\rangle$ is a strict partial order,
  $\langle P, \preceq\rangle$ is the corresponding
  non-strict partial order, and $F \subseteq P$, then $F$ is
  a filter in  $\langle P,
  \prec\rangle$ if and only if $F$ is a filter in 
 $\langle P, \preceq\rangle$, and
  \textit{vice versa}. Moreover, $F$ is a maximal
  (unbounded) filter in either of these partial orders if
  and only if it is maximal (unbounded, respectively) in the other partial
  order.
  
  A topology on the set of maximal (unbounded) filters of a
  strict partial is defined in the same way as for a
  non-strict partial order.  Once this definition is made, it
  is immediate that for any strict poset $\langle P,
  \prec\rangle$ and corresponding non-strict poset $\langle
  P, \preceq\rangle$, the identity map $P \to P$ induces a
  homeomorphism of the topological spaces of maximal
  (unbounded, respectively) filters of these posets.  For
  this reason, when it is convenient, we may prove results using
  strict partial orders instead of non-strict partial
  orders. This technique is sound because any example of
  a poset space obtained from a strict partial order can be
  converted to a homeomorphic example obtained from a
  non-strict poset space, and \textit{vice versa}.
\end{remark}

We now present two examples showing that many familiar
spaces are homeomorphic to poset spaces.

\begin{theorem}\label{s2t1}
  Every locally compact Hausdorff space is homeomorphic to
  an MF space.
\end{theorem}

\begin{proof}
  Let $X$ be a locally compact Hausdorff space and let $P$
  be the set of all nonempty precompact open subsets of $X$.  For
  $U$,$V \in P$ we put $U \preceq V$ if $U = V$ or the
  closure of $U$ is contained in $V$.  If $F$ is a filter
  and $U \in F$ then, because $U$ is precompact, 
\[
\bigcap F = \bigcap \{
  \bar{V} \mid \bar{V} \subseteq U, V \in F\}
\]
is the filtered intersection of non-empty compact sets, and
hence is non-empty and compact. Since $X$ is Hausdorff, any
two points of $X$ have open neighborhoods whose closures are
disjoint. If $F$ is a a maximal filter, then at most one of
these neighborhoods can be in $F$, which implies that
$\bigcap F$ is a singleton. Finally, the mapping $\phi\colon
\MF(P) \to X$ given by $F \mapsto \bigcap F$ has as its
inverse the mapping
\[
\phi^{-1}\colon x \mapsto \{ p \in P \mid x \in N_p \}.
\]

To prove that $\phi$ is continuous, fix $x \in \MF(P)$ and let
$U$ be any open neighborhood of $\phi(x)$ in~$X$.  Because
$X$ is locally compact, we may assume without loss of
generality that $U$ is precompact, because the precompact
sets form a basis for the topology. Thus we assume $U = N_p$ for
some $p \in P$. Now, since $\phi(x) \in U$, we have $p \in
x$, so $x \in N_p$.  Moreover, for any $F \in N_p$ in
$\MF(P)$, we have $\phi(F) = \bigcap F \in U$. This shows
$\phi$ is continuous.

To prove that $\phi^{-1}$ is continuous, let $y \in X$ be fixed,
and let $V$ be any open neighbohood of
$\phi^{-1}(y)$ in $\MF(P)$. Without loss of generality we may assume $V
= N_p$ for some $p \in P$. Now $p$ itself is some precompact
open subset $U$ of $X$, and for any $y' \in U$ we have $p \in
\phi^{-1}(y')$. Thus $\phi^{-1}(U) \subseteq V$. 
This shows $\phi^{-1}$ is continuous.
\end{proof}

\noindent
As there are non-locally-compact complete separable metric
spaces and locally compact Hausdorff nonmetrisable spaces,
the next theorem is independent of Theorem~\ref{s2t1}.  A
construction similar to that in the next theorem was used by
Lawson~\cite{Lawson-SMP} to represent complete separable
metric spaces in the context of domain theory (see
Section~\ref{sec5a}).

\begin{theorem}\label{s2t2}
  For every complete metric space $X$ there is a poset $P$
  such that $X \cong \UF(P)$ and $\UF(P) = \MF(P)$.
  Moreover, if $X$ is infinite then we may take the
  cardinality of $P$ to be that of any dense subset of $X$.
\end{theorem}

\begin{proof}
  Let $X$ be a complete metric space; we write
  $B(x,\epsilon)$ for the open metric ball of radius
  $\epsilon > 0$ around a point $x \in X$.  Let $A$ be a
  dense subset of $X$.  The poset $P$ is the set of all open
  balls $B(a,r)$ where $r$ is a positive rational number and
  $a \in A$.  For $p = B(a,r)$ and $p' = B(a',r')$ in $P$ we
  let $p \prec p'$ if and only if $d(a,a') + r < r'$.  An
  argument similar to the one in the proof of Theorem~\ref{s2t1} shows
  that any unbounded filter on $P$ has a unique point in its
  intersection. The resulting mapping $\phi\colon F \mapsto \bigcap F$
  from $\UF(P)$ to $X$ has as its inverse the mapping 
\[
x \mapsto \{ B(a,r) \mid x \in B(a,r), a \in A, r \in
\setQ^+\}.
\]
Each of these mappings can be shown to be continuous using
the same method as the proof of Theorem~\ref{s2t1}, using the fact that 
the open balls included in $P$ form a basis for~$X$.
Finally, since $X$ is a complete metric space, every
unbounded filter is maximal (see Theorem~\ref{sec3t1} below
for details).
\end{proof}

\noindent If Theorem~\ref{s2t2} is applied to the real line
using the line itself as the dense subset, the resulting
poset $P$ will be uncountable, but $\MF(P) = \UF(P)$ will be
homeomorphic to the real line.

There are also second-countable nonmetrisable Hausdorff MF
spaces. One example is the Gandy--Harrington space from
modern descriptive set theory (see~\cite{Mummert-MFspaces}).
 
\section{Separation and countability properties}\label{sec3}

\noindent
In this section, we determine the separation properties that
a poset space must satisfy.  We then show that every second
countable poset space is homeomorphic to a poset space
obtained from a countable poset.  In Section~\ref{sec6}, we
will show that a countably based poset space is either
countable or contains a perfect closed set.

\begin{theorem}\label{sec3t1}
  \begin{enumerate}
  \item Every UF space is $T_0$.
  \item Every MF space is $T_1$.
  \item If $\UF(P)$ is $T_1$ then every unbounded filter on
    $P$ is maximal and thus $\UF(P) = \MF(P)$.
  \end{enumerate}
\end{theorem}

\begin{proof}
  (1) follows from the fact that distinct filters are
  distinct as subsets of $P$. (2) follows from the fact that
  no maximal filter can properly contain another maximal
  filter.  To prove (3), suppose $\UF(P)$ is $T_1$ and let
  $F$ be an unbounded filter on $P$.  Let $G$ be a filter on
  $P$ such that $F \subseteq G$.  Clearly $G$ is unbounded.
  If $F \not = G$ then there must be a $p \in P$ such that
  $F \in N_p$ and $G \not \in N_p$.  This means $p \in (F
  \setminus G)$, which is impossible.  Thus $F = G$; this
  shows that $F$ is maximal.
\end{proof}

\begin{theorem}\label{s3t2}
  Suppose that $P$ is a poset such that $\MF(P)$ is second
  countable. There is a countable subposet $R$ of $P$ such that
  the map $F \mapsto R \cap F$ is a homeomorphism from $\MF(P)$
  to $\MF(R)$.
\end{theorem}

\begin{proof}
  Suppose that $\MF(P)$ is second countable; thus $P$
  contains a countable subset $Q_0$ such that $\{N_q: q \in
  Q_0\}$ is a basis for the topology, because every basis of
  a second-countable topology contains a countable subclass
  which is also a basis.  For $n=0,1,2,\ldots$, we construct
  a set $Q_{n+1}$ inductively to satisfy the following
  conditions.
  \begin{itemize}
  \item $Q_{n+1}$ is countable.
  \item $Q_n \subseteq Q_{n+1} \subseteq P$.
  \item For every $F \in \MF(P)$ and every finite subset $D
    \subseteq Q_n \cap F$ there is a $q \in Q_{n+1}$ such
    that $q \preceq d$ for all $d \in D$.
  \end{itemize}
  In order to see that $Q_{n+1}$ can be taken to be
  countable, suppose $D$ is a finite subset of $Q_n$ with
  nonempty intersection.  Let $E_D$
%E_D = \{p \in P:  \forall d \i D\,(p \preceq d)\}$
  be the set of all $p \in P$ such that $p \preceq d$ for
  every $d \in D$. For every filter $F \in \MF(P)$ with $D
  \subseteq F$ there is an element $p \in E_D \cap F$; thus
  $\{N_e: e \in E_D\}$ is an open cover of the intersection
  of all open sets $N_d$ with $d \in D$.  Since the given
  space is second countable, there is a countable subset
  $F_D$ of $E_D$ covering the same set of maximal filters;
  if some finite subset $D$ of $Q_n$ is not contained in any
  filter then let $F_D$ be empty.  Now take $Q_{n+1}$ to be
  the union of all $F_D$ where $D \subseteq Q_n$ and $D$ is
  finite; $Q_{n+1}$ is also at most countable.

  Let $R = \bigcup_i Q_i$. Note that $R$ is countable and
  $\{N_r \mid r \in R\}$ is a basis for $\MF(P)$.  For $F
  \subseteq P$ we write $\phi(F)$ for $F \cap R$.  It is
  straightforward to verify that $\phi(F)$ is a filter for
  every $F \in \MF(P)$, by the construction of $R$.  Because
  $R \subseteq P$, every $F \in \MF(R)$ extends to some $F'
  \in \MF(P)$; then $\phi(F') = F$.  This shows that $\phi$
  determines a surjective map $\Phi$ from $\MF(P)$ to
  $\MF(R)$.

  In order to prove that $\Phi$ is injective, it suffices to
  prove the following statement. For maximal filters $V,W$ on
  $P$ we have $V \subseteq W$ if and only if $\phi(V)
  \subseteq \phi(W)$.  Suppose $p \in V \setminus W$.  Then
  $W \notin N_p$ and thus $W \notin N_q$ for all $q$ with
  $N_q \subseteq N_p$. On the other hand, $R$ is a basis and
  $N_p$ is the union of basic open sets.  Since $V \in N_p$
  there is a $r \in R$ with $N_r \subseteq N_p$ and $V \in
  N_r$.  It follows that $r \in \phi(V)\setminus \phi(W)$.
  The other direction of the implication is trivial.

  This shows that $\phi$ is a bijection from $\MF(P)$ to
  $\MF(R)$. To see that $\phi$ is continuous, let $x \in
  \MF(P)$ be fixed and let $U$ be an open neighborhood of
  $\phi(x) = x \cap R$ in $\MF(R)$.  Without loss of
  generality we may assume that $U$ is of the form $N_r$ for
  some $r \in R$. Let $V = \{ y \in \MF(P) \mid r \in y\}$
  be the basic open set determined by $r$ in $\MF(P)$.  Now,
  because $r \in \phi(x) = x \cap R$, we see that $r \in x$,
  and thus $x \in V$.  Moreover, for any $x' \in V$, we have
  $r \in x'$, and so $r \in x' \cap R$, which means
  $\phi(x') \in U$. Thus $\phi$ is continuous.

  To see that $\phi^{-1}$ is continuous, let $V$ be any open
  subset of $\MF(P)$, and let $\phi^{-1}(y)$ be in $V$.
  Because $\{N_r \subseteq \MF(P) \mid r \in R\}$ is a
  basis for $\MF(P)$, there is some $r \in R$ with $\phi^{-1}(y) \in
  N_r \subseteq N_p$. Moreover, any $y' \in \MF(R)$ with $r
  \in y'$ will satisfy $r \in \phi^{-1}(y')$.  Thus, for $U
  = \{ y \in \MF(R) \mid r \in y\}$, we have $y \in U$ and
  $\phi^{-1}(U) \subseteq V$. This shows that $\phi^{-1}$ is
  continuous.
\end{proof}

\begin{corollary} 
  An MF space is homeomorphic to a countably based MF space
  if and only if it is second countable.
\end{corollary}

\begin{corollary}
  A UF space is homeomorphic to a countably based UF space
  if and only if it is second countable.
\end{corollary}

\begin{proof}
  Let $X = \UF(P)$ be second countable.  Construct a poset
  $R$ and a map $\phi$ in a manner analogous to the proof of
  Theorem~\ref{s3t2}.  We show that $\phi$ is a
  homeomorphism from $\UF(P)$ to $\UF(R)$.  It is clear that
  if $F \in \UF(P)$ then $\phi(F) \in \UF(R)$.  Every $G \in
  \UF(R)$ extends to some $G' \in \UF(P)$ and then $\phi(G')
  = G$.  Thus $\phi$ is well defined and surjective as a map
  from $\UF(P)$ to $\UF(R)$.  To see that $\phi$ is
  injective, suppose that $F \not = G$ are unbounded filters
  on $P$.  Without loss of generality we may assume there is
  some $p \in G \setminus F$.  There is thus some $r$ in $R
  \cap ( G \setminus F)$, because $R$ is a basis.  But $r
  \in R \cap (G \setminus F)$ implies $r \in \phi(G)
  \setminus \phi(F)$, which shows $\phi(G) \not = \phi(F)$.
  Thus $\phi$ is a bijection from $\UF(P)$ to $\UF(R)$.  The
  proof that $\phi$ is a homeomorphism is the same as in the
  proof of Theorem~\ref{s3t2}.
\end{proof}

\section{Product and subspace properties}\label{sec4}

In this section, we show that the class of MF spaces is
closed under taking $G_\delta$ subsets and arbitrary
topological products.  The class of UF spaces is closed
under taking $G_\delta$ subspaces, but it not closed under
even finite products.

\begin{theorem}  \label{s4t1}
  The class of MF spaces is closed under arbitrary
  topological products.  
\end{theorem}

\begin{proof}
  Suppose that we are given a collection $\langle \langle
  P_i, \preceq_i \rangle \mid i \in I \rangle$ of posets.
  We may assume without loss of generality that each poset
  has a greatest element, which we denote by $p_i$.  We form
  a poset $P$ consisting of those functions $f$ from $I$ to
  $\bigcup_{i \in I} P_i$ such that $f(i) \in P_i$ for all
  $i$ and $f(i) = p_i$ for all but finitely many $i$.  For
  $f,g \in P$ we put $f \preceq g$ if $f(i) \preceq_i g(i)$
  for all $i$.

  We define a map $\phi$ from $\prod_i \MF(P_i)$ to $\MF(P)$ by
  sending $\prod_i F_i$ to the set of all functions $f \in P$
  such that $f(i) \in F_i$ for all $i$. The inverse of $\phi$
  takes $x \in \MF(P)$ and returns and returns $\prod_i x_i$,
  where 
\[
x_i = \{ p \in P_i \mid \text{for some } q \in x, q(i) = p\}.
\]

To see that $\phi$ is continuous, let $x \in \prod_i
\MF(P_i)$ be fixed and let $U$ be a basic open neighborhood
of $\phi(x)$, so $U$ is of the form $N_{p}$ for some $p \in
P$.  Now $p$ is represented by a function $f \colon I \to
\bigcup_{i} P_i$ that returns the maximal element of
$P_i$ for all but finitely many $i \in I$.  Thus $f$
determines a basic open set $V$ in the product topology
$\prod_i \MF(P_i)$ such that $V$ is equal, in each coordinate~$i \in I$,
to the open set determined by $f(i)$.  Then $x \in
V$. Suppose $x' = \prod_i
x'_i$ is any point of $\prod_i \MF(P_i)$ that is in $V$,
meaning that $f(i) \in x'_i$ for all $i \in I$. Then 
$\phi(x')$ will have the property that $p_i \in x'_i$ for
each $i \in I$, which means $\phi(x') \in N_p$. Thus $\phi$
is continuous.

To see that $\phi^{-1}$ is continuous, let $y \in \MF(P)$ be
fixed, and let $V$ be any neighborhood of $\phi^{-1}(y)$ in
$\prod_i \MF(P_i)$.  By the definition of the product
topology, there is a basic open neighborhood of
$\phi^{-1}(y)$ which is obtained as a product $\prod_{i}
V_i$ of open sets $V_i \subseteq \MF(P_i)$ such that $V_i =
\MF(P_i)$ for all but finitely many $i \in I$.  Moreover, in
the finitely many coordinates where $V_i$ is a proper subset
of $\MF(P_i)$, we can find a basic open subset $N_{r(i)}
\subseteq V_i$ such that the projection of $\phi^{-1}(y)$ to
coordinate $i$ is in $N_{r(i)}$. For all $i$ where $V_i =
\MF(P_i)$ we let $r(i)$ be the greatest element of $P_i$. Now
let $f$ be the element of $P$ such that $f(i) = r(i)$ for
all $i \in I$. Then $y \in N_f$ (in $\MF(P)$), and any $y' \in
N_{f}$ will satisfy $\phi^{-1}(y) \in V$. Thus $\phi^{-1}$
is continuous.
\end{proof}

\begin{corollary}
  Every topological product of countably many countably
  based MF spaces is homeomorphic to a countably based MF
  space.
\end{corollary}

\begin{proof}
  Under these hypotheses, the poset constructed in
  Theorem~\ref{s4t1} is countable.
\end{proof}

\begin{theorem}\label{s4t2}
  The class of MF spaces is closed under taking $G_\delta$
  subspaces.
\end{theorem}

\begin{proof}
  Suppose that $\langle U_i \mid i \in \setN\rangle$ is a
  sequence of open subsets of $\MF(P)$ and $U = \bigcap_i U_i$
  is nonempty.  We form a poset $Q$ of pairs $\langle n,
  p\rangle$ such that $n \in \setN$ and $N_p \subseteq
  \bigcap_{i < n} U_i$, declaring $\langle n,p\rangle \prec
  \langle n', p'\rangle$ if $n > n'$ and $p \preceq p'$.  We
  define a map $\phi$ from $\bigcap_i U_i$ to $\MF(Q)$ by sending a
  maximal filter $F$ to the set of all $\langle n,p \rangle$
  in $Q$ such that $F \in N_p$. The inverse $\psi$ of $\phi$ takes
  a maximal filter  $G \in \MF(Q)$ and returns the set
\[
\psi(G) = \{ p \in P \mid \text{for some } n \in \setN, \langle n,p
\rangle \in x
\}
\] 
To see that $\psi(G)$ is a filter, note that if $\langle n,p
\rangle \in G$ and $\langle m,q \rangle \in G$ then there is
some common extension $\langle o, r \rangle \in G$, and thus
$r$ is a common extension of $p$ and $q$ in $\psi(G)$. 

To see that $\psi(G)$ is maximal, note that
if $\bigcap \{ p \mid p \in \psi(G)\}$ contained more than
one point of $U$, then at least one of the points of the
intersection has a basic open neighborhood $N_q$ that does
not contain one other point of the intersection. It is then
possible to adjoin $N_q$ to $G$ and extend this to a filter,
contradicting the maximality of $G$.

To see that $\phi$ is continuous, note that if $\phi(F) \in
\langle n,p\rangle$ then for every $F' \in \MF(P) \cap N_p$,
we have $\phi(F) \in \langle n,p \rangle$. Conversely, if
$\psi(G) \in N_p \cap U$ then $\psi(G) \in U_1$ and thus
every $G' \in N_{\langle 1,p\rangle} \subseteq \MF(Q)$ will have $\psi(G')
\in N_p$.
\end{proof}

\noindent
Theorem~\ref{s4t2} gives an optimal result.  We will show
below that all poset spaces have the property of Baire.  The
real line is homeomorphic to a UF space, but the $F_\sigma$
subset of rational numbers does not have the property of
Baire and therefore is not homeomorphic to a poset space.

The class of UF spaces does not enjoy the closure properties
that the class of MF spaces does.  We now give an example
showing that the class of UF spaces is not closed under
finite products.

\begin{example}
  There are two posets $P,Q$ such that $\MF(P) =
  \UF(P)$, $\MF(Q) = \UF(Q)$, but the topological product
  $\MF(P) \times \MF(Q)$ is not homeomorphic to any UF
  space.
\end{example}

\begin{proof}
  Let $\omega$ denote the least infinite countable ordinal
  and let $\omega_1$ denote the least uncountable ordinal.
  We define $P$ to be the set of functions from finite
  initial segments of $\omega$ to $\{0,1\}$ and define $Q$
  to be the set of functions from countable initial segments
  of $\omega_1$ to $\{0,1\}$.  For both posets the relation
  $\preceq$ is given by extension: $p \preceq q$ if, for all
  $\alpha$ in the domain of $q$, $p(\alpha)$ is defined and
  takes the value $q(\alpha)$.   

  We first show that $\MF(P) = \UF(P)$ and $\MF(Q)$ =
  $\UF(Q)$.  Assume that $F$ is an unbounded filter on $P$
  (the argument for $Q$ is parallel).  Then all functions in
  $F$ are compatible, that is, they do not contradict each
  other on any value in the intersection of their domains.
  There is thus a total limit function $f$, because
  otherwise there would be a first ordinal $\alpha$ where
  $f$ is undefined and the function extending $f$ which maps
  $\alpha$ to $0$ would define an element of $P$ which
  would be a lower bound for the filter $F$. Since $f$ is total, all
  functions mapping the ordinals up to some $\alpha$ in the
  domain of $f$ to the corresponding value of $f$ are in the
  filter.  One can see that this filter is already maximal,
  because any element outside it but still in $P$ is
  incompatible with this function and adding it would
  destroy the filter property.

  Assume now, by way of contradiction, that $\UF(P) \times
  \UF(Q)$ is homeomorphic to a space $\UF(R)$. We denote by
  $\pi_P, \pi_Q$ the continuous, open projection maps from
  $\UF(R)$ to its factor spaces.  There is a filter $F$ in
  $\UF(R)$ such that $\pi_P(F)$ and $\pi_Q(F)$ are the
  filters generated by the set of all functions in $P$ and
  $Q$, respectively, which map all inputs to $0$.  Now one
  can select an infinite sequence $r_0,r_1,\ldots$ in $F$
  such that for each $n$ the projection $\pi_P(N_{r_n})$
  consists only of functions which map the first $n$ numbers
  to $0$ and $r_{n+1} \preceq r_n$ for all~$n$. The sequence
  $\langle r_i\rangle$ generates a subfilter $G \subseteq F$.
  There is no lower bound~$r$ for $G$, because otherwise
  $\pi_P(N_r)$ would be an open set containing some basic
  open set $N_p$ such that $N_p \subseteq \pi_P(N_{r_n})$ for
  all~$n$; such a $p$ cannot exist by construction.

  On the other hand, there is a function $f$ contained in
  all the open sets $\pi_Q(N_{r_n})$ and there are basic
  open neighbourhoods of $f$ generated by $q_0,q_1,\ldots$
  such that $N_{q_n} \subseteq \pi_Q(N_{r_n})$ for each $n$.
  The basic open sets $N_{q_0},N_{q_1},\ldots$ fix $f$ only on
  countably many ordinals and thus their intersection is
  also a basic open set.  So $\pi_Q(G)$ is bounded while
  $\pi_Q(F)$ is not and thus $G \subset F$.  It follows that
  $\UF(R)$ is not a $T_1$ space.  This contradicts the assumption 
  that $\UF(R)$ is homeomorphic to $\MF(P) \times \MF(Q)$.
\end{proof}

\noindent
We note that the previous example is not second countable
and that the failure of second countability was important to
the proof.

\begin{question}
  Is the class of countably based UF spaces closed under
  taking finite (or arbitrary) topological products?
\end{question}

\noindent
We end the section by showing that the class of UF spaces is
closed under taking $G_\delta$ subspaces. As with the class
of MF spaces, this result cannot be extended to include
$F_\sigma$ subspaces. We first prove the result for open
subspaces, which has a much simpler proof.

\begin{theorem}
  The class of UF spaces is closed under taking open subspaces.
\end{theorem}

\begin{proof}
  Let $P$ be a poset and let $U$ be an open subset of $\UF(P)$. Let
  $R$ be the set of all $r \in P$ such that $N_r \subseteq U$; we
  regard $R$ as a subposet of $P$.  Then any $x \in U$ has a
  neighborhood $N_r \subseteq U$, where $r \in R$.  Thus the
  restriction map $\phi \colon x \mapsto x \cap R$ sends each element
  of $U$ to a filter on $R$. Note that if this filter were not
  unbounded as a subset of $R$ then it has a lower bound in $R$ and
  consequently would not be unbounded in $P$.

  The inverse map of $\phi$ sends each maximal filter on $R$
  to it upward closure in $P$. If $\phi(G)$ were bounded
  below by $p \in P$, then in particular $p \preceq r$ for
  some $r \in R$. Thus $N_p \subseteq N_r \subseteq
  R$, which means $p \in R$ and $G$ is not in $\UF(R)$.
  
  To set that $\phi$ and its inverse are
  continuous, note that 
\[
\left \{ N_r^P = \{ F \in \MF(P) \mid r \in
  F\} 
\:\mid r \in R
\right \}
\]
 is a basis for the restriction of
  $\MF(P)$ to the subspace $U$, that
\[
\left \{ N_r^R = \{ F \in
  \MF(R) \mid r \in F\} 
\: \mid r \in R
\right \} 
\]
is a basis for
  $\MF(R)$, and that a point $x \in U \subseteq \MF(P)$ is in
  $N^P_r$ if and only if $\phi(x)$ is in $N_r^R$.
\end{proof}

\begin{theorem}
  The class of UF spaces is closed under taking $G_\delta$ subspaces.
\end{theorem}

\begin{proof}
Let $G_0$ be the space $\UF(P)$ for some poset $P$ with order $\prec_P$
and let $G$ be a $G_{\delta}$ subset of $G_0$. Thus there is
a descending sequence $G_1,G_2,\ldots$ of open subsets of $G_0$
such that $G_0 \supseteq G_1 \supseteq G_2 \supseteq \ldots$ and
$G = \bigcap_n G_n$. Define
\[
   R = \{p \in P \mid p \in F \text{ for some } F \in \UF(P)
   \cap G\}.
\]
For each $p \in R$, let $g(p) = \sup \,\{n \in \setN\mid N_p \subseteq
G_n\}$, where $g(p)=\infty$ if $N_p \subseteq G$.
Define an order relation $\prec_R$ on
$R$ by putting $p \prec_R q$ if $p \prec_P q$ and either $g(q)
< g(p) \leq \infty $ or $g(q) = g(p) = \infty$. We will show that the
unbounded filters on $(R,\preceq_R)$ are precisely the
unbounded filters on $P$ that are in $G$ and do this by showing
the following four claims.

\textit{Claim 1:} Let $F \in G \subseteq \UF(P)$; then $F$ is an
unbounded filter in $R$ under $\preceq_R$.
By definition of $R$, $F \subseteq R$.
To show that $F$ is a filter on $R$, fix $p,q \in F$. If
$g(p)$ or $g(q)$ is infinite then $p$ and $q$ have a common
extension $r$ under $\preceq_P$ with $g(r) = \infty$. Thus
$r$ is a common extension of $p$ and $q$ under $\preceq_R$.
Otherwise, because $F \in G$, there is an $r \in F$ with
$r \preceq_P p$, $r \preceq_P q$ and $N_r \subseteq G_{g(p)+g(q)+1}$.
Then $g(r)>g(p)+g(q)$, $r \prec_R p$ and $r \prec_R q$.
As $\preceq_R$ is a restriction of $\preceq_P$, $F$ is
upward closed under $\preceq_R$ and $F$ is a filter in $R$.
Furthermore, $F$ must be unbounded in $R$, because a bound in $R$
would also be in a bound $P$.

\textit{Claim 2:} Let $F \subseteq R$ be a filter in $R$;
then either $\sup \,\{ g(p) \mid p \in F\} = \infty$ or $F$
is bounded. Suppose the
supremum is $n < \infty$ instead. There can only be one $r
\in F$ with $g(r) = n$, because $F$ is a filter on $R$. Because $r
\in R$, there is some $F' \in \UF(P)$ with $r \in F'$ and
$F' \in G$. Thus there is an $r' \in F'$ with $g(r') >
g(r)$ and $r' \prec_P r$; this means $r' \prec_R r$,
which shows that $F$ is bounded in $R$. 
%This contradiction establishes the claim.

\textit{Claim 3:} Let $F$ be a bounded filter of $P$ which is
also a filter in $R$; then $F$ is bounded in $R$. Let $r \in P$
be a lower bound for $F$. If $N_r \not\subseteq G_n$ for
some $n$ then $\sup \,\{ g(p) \mid p \in F\} < n$ and $F$ is
bounded in $R$ by Claim 2. Otherwise $N_r \subseteq G$,
in which case $r \in R$ and $F$ is again bounded as a subset of $R$.

\textit{Claim 4:} Let $F$ be an unbounded filter in $R$; then
$F$ is also an unbounded filter in $P$.
To see this, consider the upward closure $F'$ of $F$ in $P$.
$F'$ is unbounded in $P$, by Claim~3.
Claim~2 shows that $F' \in G$; thus $F' \subseteq R$. The
definition of $F'$ shows that $F \subseteq F'$. Fix $r \in F'$;
then there must be a $p \in F$ with $p \preceq_P r$. If
$g(p) = \infty$ then $p \preceq_R r$ and so $r \in F$.
Otherwise there must be a $q \in F$ with $q \preceq_P p$ and
$g(q) > g(r)$. Then it follows from transitivity of $\preceq_P$
and the definition of $\prec_R$ that $q \prec_P r$, $q \prec_R r$
and $r \in F$. This shows $F' = F$.

Claims 1 and 4 show that the unbounded filters on $R$ are
exactly those unbounded filters on $P$ which are in $G$.  So
the identity map $\phi \colon \UF(P) \cap G \to \UF(R)$ is
surjective by Claim~4. As this map is trivially injective, it is
thus invertible.
To see that $\phi$ and $\phi^{-1}$ are continuous, let $x
\in \UF(P) \cap G$ be fixed.  Note that for any $r \in R$,
we have $r \in x$ if and only if $r \in \phi(x)$, because
$\phi$ is the identity map on filters. Thus $\phi(x)$ is in
the basic open neighorhood of $\UF(R)$ determined by $r$ if
and only if $x$ is in the basic open neighborhood of $\UF(P)
\cap G$ determined by $r$. 
\end{proof}

\section{Completeness properties}\label{sec5}

\noindent
In this section, we establish that every poset space has the
a completeness property known as the strong Choquet
property. We then characterize the class of countably based MF spaces
as precisely the class of second-countable $T_1$ spaces with
the strong Choquet property. We first establish a weaker
property.

\begin{theorem}  
  Every poset space has the property of Baire.
\end{theorem}

\begin{proof}
  Let $X$ be $\MF(P)$ or $\UF(P)$. Suppose that $\langle
  U_i \mid i \in \setN \rangle$ is a sequence of dense open
  sets in $X$ and $V$ is a fixed open set. We construct a
  sequence $\langle p_i \mid i \in \setN\rangle$ of elements 
  of $P$. Let
  $p_0$ be such that $N_{p_0} \subseteq V \cap U_0$. Given
  $p_i$, there is an unbounded or maximal filter in $N_{p_i}
  \cap U_{i+1}$. Choose $p_{i+1}$ such that $N_{p_{i+1}}
  \subseteq U_{i+1} \cap N_{p_i}$ and $p_{i+1} \preceq p_i$.
  In the end, $F = \langle p_i \rangle$ is a linearly
  ordered subset of $P$. Thus $F$ extends to an element of
  $X$. Clearly this element is in $V \cap \bigcap_i U_i$.
\end{proof}

\noindent
We will now show that every poset space has the strong
Choquet property, which is defined using a certain game
first introduced by Choquet~\cite{C-LA}. Let $X$ be an arbitrary
topological space. The strong Choquet game is the
Gale--Stewart game (see \cite{GS-IGPI}
and~\cite{Kechris-CDST}) defined as follows.  The stages of
play are numbered $0,1,2,\ldots$ and both players make a
move in each stage. In stage $i$, player~\I plays an open
set $U_i$ and a point $x_i$ such that $x_i \in U_i$ and if
$i>0$ then $U_i \subseteq V_{i-1}$.  Then player~\II plays
an open set $V_i$ such that $x_i \in V_i$ and $V_i \subseteq
U_i$. At the end of the game, \iffalse At stage $0$,
player~\I plays an open set $U_0$ and a point $x_0 \in U_0$.
At stage 0, player~II plays an open set $V_0$ such that $x_0
\in V_0 \subseteq U_0$. At stage $i+1$, player~\I plays an
open set $U_{i+1}$ and a point $x_{i+1}$ such that $x_{i+1}
\in U_{i+1} \subseteq V_{i}$. At stage $i+1$, player~II
plays an open set $V_{i+1}$ such that $x_{i+1} \in V_{i+1}
\subseteq U_{i+1}$. At the end of the game, \fi player~\I
wins if $\bigcap_i U_i$ is empty (or, equivalently, if
$\bigcap_i V_i$ is empty). Player~\II wins if $\bigcap_i
U_i$ is nonempty. A \textit{position} in the game is a
finite (possibly empty) sequence
\[
\langle \langle U_0, x_0 \rangle, V_0, \langle U_1, x_1
\rangle, \ldots \rangle
\]
which is an initial segment of an infinite play of the game
following the rules just described.

A space $X$ has the \define{strong Choquet property} if
player~II has a winning strategy for the strong Choquet game
on $X$. A winning strategy is a function that takes a
position after player~\I has played and tells player~II which
open set to play, such that if player~II follows the winning
strategy then player~II will always win the game regardless
of what moves player~\I makes.

The strong Choquet property is strictly stronger than the
property of Baire. Moreover, the class of topological
spaces with the strong Choquet property is closed under
$G_\delta$ subspaces and arbitrary topological products. It
is known that the class of topological spaces with the
property of Baire is not closed under binary products
(an example is provided in~\cite{FK-BBS}).

\begin{theorem}\label{s5t1}
  Every poset space has the strong Choquet property.
\end{theorem}

\begin{proof}
  We describe the strategy for player~II informally. At the
  start of the game, player~\I plays an open set $U_0$ and a
  point $x_0$. Player~II translates the point $x_0$ into a
  filter on $P$, then finds a basic neighbourhood $q_0$ of
  $x$ such that $N_{q_0} \subseteq U_0$. Player~II then
  plays $N_{q_0}$. Now given $\langle x_1, U_1\rangle$ with
  $x_1 \in N_{q_0}$, Player~II translates $x_1$ to a filter on $P$
  and then finds a neighbourhood $q_1$ of $x_1$ such that
  $q_1 \preceq_P q_0$ and $N_{q_1} \subseteq U_1$. Player
  II plays $N_{q_1}$. Player~II continues this strategy,
  always choosing $q_{i+1} \preceq_P q_i$. At the end of
  the game, player~II has determined $\{q_i \mid i \in \setN\}$, a
  descending sequence of elements of $P$. This sequence 
  extends to an element of $X$ which is in $\bigcap
  N_{q_i}$. Player~II has thus won the game.
\end{proof}

\noindent
We use the strong Choquet property to obtain the following
characterization of countably based MF spaces. 

\begin{theorem}\label{s5t2}
  A topological space is homeomorphic to a countably based
  MF space if and only if it is second countable, $T_1$ and
  has the strong Choquet property.
\end{theorem}

\noindent
We postpone the proof of this theorem temporarily to comment
on the hypotheses involved in the characterization.
Clearly, any space $X$ homeomorphic to a countably based MF
space must be $T_1$ and second countable. We have already
shown $X$ must also have the strong Choquet property. Thus
the new content of Theorem~\ref{s5t2} is that the strong
Choquet property is sufficient for a $T_1$ second-countable
space to be homeomorphic to a countably based MF space.
In the non-second-countable setting, the strong Choquet
property is not sufficient for a $T_1$ space to be
homeomorphic to an MF space.  

\begin{example}\label{s5e12}
  There is a Hausdorff strong Choquet space which is not
  homeomorphic to any MF space.
\end{example}

\begin{proof}
  The space $X$ consists of certain functions from
  $\omega_1$ to $\{0,1\}$. We put a function $f$ in $X$ if
  and only if there is an ordinal $\alpha < \omega_1$ such
  that $f(\beta) = 0$ for all $\beta >\alpha$. For each $f
  \in X$ and each $\alpha < \omega_1$, the set 
\[
\{ g \in X
  \mid f(\beta) = g(\beta) \text{ for all } \beta <
  \alpha\}
\]
 is declared to be an open set. The topology
  on $X$ is the one generated by these open sets. It is
  clear that $X$ is a Hausdorff space.

  It is easy to show that $X$ has the strong Choquet
  property, as follows. All that player~II has to do is to play any basic
  open subset of the open set played by player~\I which also contains
  the point given by player~\I. In the end, the open
  sets played by player~\I in the countable number of rounds
  of the game and each round will fix countably many coordinates
  of a function in $X$. In the limit, countably many coordinates
  are fixed and we can find a point in the intersection of
  the sets played by~\I by forcing the remaining coordinates
  to map to $0$.

  We now show that $X$ is not homeomorphic to any MF space.
  Suppose, by way of contradiction, that $X \cong \MF(P)$.
  We construct a transfinite sequence $\langle p_\alpha \mid \alpha <
  \omega_1\rangle$ inductively. Let $p_0$ be
  any basic open neighbourhood of the constant $0$ function.
  Given $\langle p_\alpha \mid \alpha < \beta \rangle$,
  there is a first coordinate $\gamma < \omega_1$ which is
  not fixed by any $p_\alpha$; let $f$ be the function which is
  $0$ except at $\gamma$, and $f(\gamma) = 1$. Note that any intersection
  of countably many open sets in $X$ is open. Thus we may
  choose $p_\beta \in P$ such that $p_\beta \preceq
  p_\alpha$ for all $\alpha < \beta$ and $f \in
  N_{p_{\beta}}$. Choose any such $p_\beta$. At the end of
  this construction, $\langle p_\alpha \mid \alpha <
  \omega_1\rangle$ is linearly ordered and thus extends to a
  maximal filter $F$. Now the element of $X$ corresponding
  to $F$ sends uncountably many ordinals to $1$, which is
  impossible.
\end{proof}

\noindent
We now return to the proof of Theorem~\ref{s5t2}, which will
occupy the remainder of this section. Let $X$ be a fixed
$T_1$ space with a fixed countable basis and a fixed winning
strategy for player~II in the strong Choquet game. Our
first step is to define a poset $P$. The elements of $P$ are
called \textit{conditions}. A condition is a finite list of the form
\[
\langle A, \pi_1, \pi_2, \ldots, \pi_k \rangle
\]
satisfying the following requirements.
\begin{enumerate}
\item The set $A$ is a nonempty basic open set from the fixed
  countable basis. For each condition $c$ we let $S(c)$
  denote the basic open set $A$ appearing in~$c$.
\item Each $\pi_i$ is a finite (that is, partial) play of
  the strong Choquet game on $X$ following the fixed winning
  strategy $s_{\text{II}}$ for player~II. We require each
  $\pi_i$ to be of the form
  \[
  \begin{split}
    \langle V_1, x_1, s_{\text{II}}(V_1,x_1), V_2,
    x_2, s_{\text{II}}(V_1,x_1,V_2,x_2), \ldots,\\
    V_r, x_r, s_{\text{II}}(V_1,x_1,V_2,x_2, \ldots, V_r,
    x_r) \rangle.
  \end{split}
  \]
  Thus each $\pi_i$ ends with an open set, which we will
  denote by $U(\pi_i)$. It is allowable that $\pi$ is the
  empty sequence $\langle \rangle$, in which case $U(\pi) =
  X$.

\item If a play $\pi$ is in a condition then so is every
  initial segment of $\pi$ that ends with a move by player
  II.
\item $A \subseteq U(\pi_i)$ for each $i \leq k$. 
\end{enumerate}
We define the order $\prec$ on $P$ as follows. Let $ c =
\langle A, \pi_1, \pi_2, \ldots, \pi_k \rangle$ and $c' =
\langle A', \pi'_1, \pi'_2, \ldots, \pi'_l \rangle$ be any
two conditions. We let $c' \prec c$ if and only if
\begin{enumerate}
\item[(5)] For each finite play $\pi_i$ in $c$ there is a point
  $x_n \in S(c)$ such that the longer play
  \[
  \pi_i \smallfrown \langle A, x_n , s_{\text{II}}(\pi_i
  \smallfrown \langle A,x_n\rangle )\rangle
  \]
  is in $c'$, that is, equals $\pi'_j$ for some $j \leq l$.

\item[(6)] $A' \subseteq A$ (this is actually a consequence of
  requirement (5)).
\end{enumerate}
\noindent
Requirement~(3) in the definition of a condition allows us
to prove that the order on $P$ is transitive. Because each
condition is finite, requirement~(5) in the definition of
the order relation ensures $c \not \prec c$ for all $c \in
P$. Thus $\prec$ is a partial order on $P$.

\begin{lemma}\label{lem1}
  For any filter $F$ on $P$ the intersection $\bigcap_{c \in
    F} S(c)$ is nonempty.
\end{lemma}

\begin{proof}
  Let $\langle A_i \mid i \in \setN \rangle$ be an
  enumeration of all of the basic open sets which appear as
  $S(c)$ for some $c \in F$; here we are using the fact that
  $X$ is second countable and that each $S(c)$ is drawn from
  a fixed countable basis of $X$. It is immediate that
  $\bigcap_{c \in F} S(c)$ equals $\bigcap_{i \in \setN}
  A_i$. We will show the latter intersection is nonempty.

  We inductively construct a descending sequence of
  conditions $\langle c_i \mid i \in \setN\rangle$ and a sequence of finite
  plays $\langle \pi_i \mid i \in\setN\rangle$ so that $\pi_{i+1}$ is an
  immediate extension of $\pi_i$ for each $i \in \setN$. At
  stage $0$ let $c_0$ be any condition in $F$ such that
  $S(c_0) = A_0$ and let $\pi_0$ be any finite play in
  $c_0$.

  At stage $i + 1$ let $c$ be any condition in $F$ such
  that $S(c) = A_i$. Let $c_{i+1}$ be a common extension of
  $c$ and $c_i$ in $F$. It is clear that $S(c_{i+1})
  \subseteq S(c) = A_i$. Choose $\pi_{i+1}$ to be any play
  in $c_{i+1}$ which is an immediate extension of $\pi_i$.

  Now assume the entire sequence $\langle \pi_i\rangle$ has
  been constructed. These partial plays determine an infinite play 
  $\gamma$ of the strong Choquet game following the strategy for player~\II.  
  Thus the intersection of the open sets played by player~\I in $\gamma$ 
  is nonempty. By construction, each set $A_i$ has a subset played by 
  player~\I at some stage of $\gamma$. Thus $\bigcap_i A_i$ is nonempty.
\end{proof}

\begin{lemma}\label{lem2}
  Let $c_1$ and $c_2$ be two conditions and let $x \in
  S(c_1) \cap S(c_2)$. There is a condition $c$ such that
  $c \prec c_1$, $c \prec c_2$ and $x \in S(c)$.
\end{lemma}

\begin{proof}
  Begin by letting $c$ be empty. For each $\pi$ in $c_1$ we
  put the longer play
  \[
  \pi \smallfrown \langle S(c_1), x, s_{\text{II}}(\pi
  \smallfrown \langle S(c_1), x\rangle )\rangle
  \]
  into $c$. For each $\pi$ in $c_2$ we put
  \[
  \pi \smallfrown \langle S(c_2), x, s_{\text{II}}(\pi
  \smallfrown \langle S(c_2), x\rangle )\rangle
  \]
  into $c$.
  For each $\pi$ that has been added to $c$ we add all
  initial segments of $\pi$ ending with a move by player~\II.
  We then let $S(c)$ be a basic open neighbourhood of $x$ which is
  a subset of the open set $\bigcap_{\pi \in c} U(\pi)$.
  This construction ensures that $c$ is a condition
  satisfying the conclusions of the lemma.
\end{proof}

 \begin{lemma} \label{lem3} Let $F$ be a maximal filter on
   $P$. The intersection $\bigcap_{c \in F} S(c)$ contains a
   single point.
 \end{lemma}
 
\begin{proof}
   By Lemma~\ref{lem1} we know that $\bigcap_{c \in F} S(c)$
   is nonempty. Suppose that $x,y$ are distinct points in
   $\bigcap_{c \in F} S(c)$. Let $A$ be a basic open
   neighbourhood of $x$ such that $y \not \in A$. We
   construct a filter $G$ inductively. At stage $n$ we
   construct $G_n \subseteq P$ and in the end we let $G$ be
   the upward closure of $\bigcup_n G_n$. To begin, let
   $G_0$ = $F \cup \{ \langle A, \langle \rangle \rangle \}$. At
   stage $i + 1$, we know by induction that $x \in S(c)$ for
   every $c \in G_i$. Thus we can apply Lemma~\ref{lem2}
   repeatedly so that $G_i \subseteq G_{i+1}$, every pair of
   conditions in $G_i$ has a common extension in $G_{i+1}$
   and $x \in S(c)$ for every $c \in G_{i+1}$.

   It is immediate from the construction that $G = \bigcup_i
   G_i$ is a filter which properly extends $F$. This shows
   that $F$ was not maximal.
 \end{proof}

\begin{proof}[{\upshape \bfseries Proof of Theorem~\ref{s5t2}}]
  For each $F \in \MF(P)$ we denote the single point in
  $\bigcap_{c \in F} S(c)$ by $\phi(F)$. We show that
  $\phi$ is a homeomorphism from $\MF(P)$ to~$X$.

  We first show that $\phi$ is an injective map. Suppose
  that $F$ and $F'$ are maximal filters on $P$ such that $x \in
  \bigcap_{c \in F} S(c)$ and $x \in \bigcap_{c \in F'}
  S(c)$. By following a procedure similar to the proof of
  Lemma~\ref{lem3} we may find a filter $G$ such that $F
  \subseteq G$ and $F' \subseteq G$. Thus, by maximality,
  we have $F = F' = G$.

  We next show that $\phi$ is a surjective map. Let $x \in
  X$ be fixed. Let $\langle A_i \mid i \in \setN\rangle$ be
  a sequence of basic open sets such that $\bigcap_i A_i =
  \{x\}$. The existence of this sequence requires that $X$
  be $T_1$ and first countable. For each $i \in \setN$ let
  $c_i = \langle A_i, \langle \rangle \rangle$. Following a
  method similar to the proof of Lemma~\ref{lem3}, we can
  construct a filter $F$ such that $c_i \in F$ for each $i
  \in \setN$. Let $G$ be an extension of $F$ to a maximal
  filter. Now $S = \bigcap_{c \in G} S(c)$ is nonempty by
  Lemma~\ref{lem1} and $S \subseteq \bigcap_i A_i = \{x\}$
  by construction, so $\phi(G) = x$.

  It remains to show that $\phi$ is open and
  continuous. This follows from Lemma~\ref{lem2}; for each
  $x \in X$ and each condition $c$, we have
$c \in \phi^{-1}(x)$ if and only if  $x \in S(c)$.
  This shows that $X$ is homeomorphic to $\MF(P)$. By 
  Theorem~\ref{s3t2}, we may find a countable subposet $R$ 
  of $P$ such that $X$ is homeomorphic to $\MF(R)$.
  This completes the proof.
\end{proof}

\section{An application to domain theory}\label{sec5a}

\noindent
In this section, we apply the characterization of countably
based MF spaces to characterize those second-countable
spaces with a domain representation. Our result gives a
complete solution to the so-called model problem for
second-countable spaces in domain theory.

A domain is a certain type of poset (defined below) and
every domain is a topological space with a topology known as
the Scott topology. A domain representation of a
topological space $X$ is a domain $D$ such that $X$ is
homeomorphic to the topological space consisting of the
maximal elements of $D$ with the relative Scott topology.
The history of such representations is thoroughly described
by Martin~\cite{Martin-TGDT}. It is known that every
complete separable metric space has a domain representation
(see Lawson~\cite{Lawson-SMP}) and that every space with a
domain representation is $T_1$ and has the strong Choquet property
(Martin~\cite{Martin-TGDT}). We now show that the strong
Choquet property is sufficient for a $T_1$ second-countable
space to have a domain representation.

\newcommand{\oDownarrow}{\mathord{\Downarrow}}
\newcommand{\oUparrow}{\mathord{\Uparrow}}

We summarize the definitions from domain theory that we
require; these definitions are explored fully by Gierz
\textit{et al.}~\cite{GHKLMS-CLD}. A nonempty subset $I$ of a poset
$\langle P, \preceq\rangle$ is \textit{directed} if every
pair of elements in~$I$ has an upper bound in~$I$. A poset
$P$ is said to be a \textit{dcpo} (for ``directed-complete
partial ordering'') if every directed subset of
$P$ has a least upper bound. Any dcpo $D$ has a second
order relation $\ll$, known as the \textit{way below} relation,
under which $q \ll p$ if and only if
whenever $I \subseteq D$ is a directed set with $p \ll \sup
I$ there is some $r \in I$ with $q \preceq r$. For each $p
\in D$ we put $\oDownarrow p = \{ q \in D \mid q \ll p\}$
and $\oUparrow q = \{ p \in D \mid q \ll p\}$. A dcpo
$D$ is \textit{continuous} if $\oDownarrow p$ is directed and
the equality $p = \sup
\oDownarrow p$ holds for every $p \in D$. A \textit{domain}
is a continuous dcpo. A subset $B$ of a domain $D$ is a
\textit{basis} if $B \cap \oDownarrow p$ is directed and $p
= \sup (B \cap \oDownarrow p)$ for every $p \in D$. A domain
is \mbox{\textit{$\omega$-continuous}} if it has a countable
basis. An element $p$ of a dcpo is \textit{compact} if
$p \ll p$. A dcpo $D$ is \textit{$\omega$-algebraic}
if there is a countable basis for $D$ consisting of compact
elements.
The \textit{Scott topology} on a dcpo $D$ is
generated by the basis $\{ \oUparrow p \mid p \in D\}$. A
\textit{domain representation} of a space $X$ is a
homeomorphism between $X$ and the maximal elements of a
domain with the Scott topology.

\begin{theorem}
  A topological space has a domain representation via
  an \hbox{$\omega$-algebraic} dcpo if and only if the space
  is second-countable, $T_1$ and has the strong Choquet property.
\end{theorem}

\begin{proof}
  It can be seen that any space with a domain representation
  satisfies the $T_1$ separation property and a result of
  Martin~\cite{Martin-TGDT} shows that any space with a
  domain representation has the strong Choquet property.
  Therefore, we only need to prove that a second-countable
  $T_1$ strong Choquet space has a domain representation via
  an $\omega$-algebraic dcpo We use the following
  lemma, which follows easily from the definitions.

\begin{lemma}\label{s6l12}
  Suppose that $P$ is a countable poset. The set of all
  filters on~$P$, ordered by inclusion, is an
  $\omega$-algebraic dcpo $D$. The maximal filters on
  $P$ are precisely the maximal elements of $D$ and the
  compact elements of $D$ are precisely the principal
  filters on $P$. Moreover,
  the poset topology on $\MF(P)$ corresponds exactly to the
  Scott topology on the maximal elements of $D$.
\end{lemma}

\noindent
We showed in Theorem~\ref{s5t2} that any second-countable
$T_1$ strong Choquet space is homeomorphic to $\MF(P)$ for a
countable poset $P$. It follows immediately from the lemma
above that such a space also has a domain representation
via an $\omega$-algebraic dcpo.
\end{proof}

\noindent The next corollary follows from the fact that any space with a
domain representation is $T_1$ and has the strong Choquet property. 
Although this corollary is already known, the proof here is new.

\begin{corollary}
  If a second-countable space has a domain representation
  then it has a representation via an $\omega$-algebraic dcpo.
\end{corollary}

\noindent We end this section with several remarks on the relationship
between domain representable spaces and MF spaces.

A proof of Lemma~\ref{s6l12} can be modified to show that the collection of
all ideals on a poset (sometimes called the \textit{ideal completion} of
the poset) forms a domain whose maximal elements in the Scott topology
correspond to the maximal ideals of the poset under the Stone topology. 
All results we have proved for MF spaces also hold for these spaces of
maximal ideals, by duality.  The relationship between ideal completions and 
domain representations has been investigated by Martin~\cite{Martin-IMS}.

A \textit{Scott domain} is a domain in which every pair of
elements with an upper bound has a least upper bound.
Lawson~\cite{Lawson-SMP} has shown that any space with a
domain representation via a countably based Scott domain is
a complete separable metric space. It can be seen that
posets constructed in Theorem~\ref{s5t2} do not, in general,
give Scott domains, even when the posets are constructed
from formal balls in complete separable metric spaces.

The proof of Example~\ref{s5e12} can be modifed to obtain 
the following.
\begin{example}
There is a Hausdorff strong Choquet space that does not have a domain 
representation.
\end{example}

\section{Semi-Topogenous Orders}\label{sec7}

\noindent In this section, we prove results which give a partial
solution to the question of which arbitrary (not necessarily
second countable) topological spaces are homeomorphic to MF spaces.

Suppose that a topological space $X$ is homeomorphic to
$\MF(P)$, for some poset $P$, via a fixed homeomorphism
$\phi$.  If each element of $p \in P$ is replaced by the
corresponding open subset $\phi(N_p) \subseteq X$, the poset
order on $P$ will determine a corresponding order relation
on these subsets of~$X$. Moreover, the collection of all
these open subsets forms a basis for the topology on $X$. It
is thus natural to ask whether the existence of a basis with
a suitable order relation is sufficient for
a topological space to be homeomorphic to an MF space.

Cs\'asz\'ar~\cite{C-FGT} considered many different types of
orders and their connections to topology. The basic concept
is that of a semi-topogenous order.

\begin{definition}
A \textit{semi-topogenous order} is a binary relation 
$\sqsubset$ on the powerset of a topological space $X$ satisfying the 
following axioms for all $u,v,w \subseteq X$ \cite[Chapter 2]{C-FGT}:
\begin{itemize}
\item $\varnothing \sqsubset \varnothing$ and $X \sqsubset X$;
\item $v \sqsubset w \Rightarrow v \subseteq w$;
\item $u \subseteq v \sqsubset w \Rightarrow u \sqsubset w$;
\item $u \sqsubset v \subseteq w \Rightarrow u \sqsubset w$.
\end{itemize}
\end{definition}

\noindent
Cs\'asz\'ar considered orders which are only linked to
topology, such as the order which says that $w$ is a
neighbourhood of $v$.
%, that is, that $v$ is not near to $X \setminus w$.
It might happen that some but not all open
supersets $w$ of a given set $v$ satisfy $v \sqsubset w$.
Nevertheless, although this is not made explicit by
Cs\'asz\'ar, it is quite convenient to postulate also a
connection between the topology and the open spaces.

Recall that the \textit{open kernel} of a set is the union of all its
open subsets. We say that the topological space $X$ is
\textit{generated} by the order $\sqsubset$ if for each $u \subseteq X$ the
set $\bigcup \{ o \subseteq X \mid o \sqsubset u \}$ is the open kernel of $u$.
In this case a set $w$ is open if and only if it is the union of all
$v$ such that $v \sqsubset w$. It follows that if $v \sqsubset w$ then
there is an open $o$ with $v \subseteq o \subseteq w$; the converse of this
last implication does not always hold. Every topological
space is generated by some semi-topogenous order, for one
can define $v \sqsubset w$ to hold if and only if there is an open
set $o$ with $v \subseteq o \subseteq w$.

\begin{remark} 
  There is a close relationship between semi-topogenous
  orders and the way below relation $\ll$ on a continuous
  dcpo, which was discussed in Section~\ref{sec5a}. The
  following true properties of the way below relation are
  obtained by dualizing the second, third and fourth
  properties in the definition of a semi-topogenous order:
\begin{align*}
v \ll w \Rightarrow v \leq w\\
u \leq v \ll w \Rightarrow u \ll w\\
u \ll v \leq w \Rightarrow u \ll w
\end{align*}
The fact that these are dual forms follows from the fact
that points in a topological space are minimal as nonempty subsets
under $\subseteq$ but are maximal elements of a domain representing
the topological space; for this reason, we write $\leq$ for $\supseteq$
and $\ll$ for $\sqsupset$.
The requirement that $\bigcup \{o \mid o \sqsubset u\}$ is
the open kernel of $u$ corresponds exactly to the fact that
$\{x \mid y \ll x\}$ is the open kernel of an element $y$ of a
continuous dcpo with the Scott topology.

Thus if a space $X$ has a representation via a continuous
dcpo $D$ then the dual of the way below relation on $D$
is a semi-topogenous order (except that it is defined only
on a subset of the powerset of $X$) which generates the
topology on $X$. Semi-topogenous orders can be viewed as a
generalization of the way below relation which is applicable
to the case when the dcpo is the full powerset of a
topological space. It appears that semi-topogenous orders are
related to auxilliary relations as defined by Gierz 
\textit{et al.}~\cite{GHKLMS-CLD}, although a formal relationship seems 
difficult 
to state.
\end{remark}

\noindent
A \textit{filter} in a topological space $X$ is a collection
of nonempty subsets that is closed under finite intersection
and under superset. A filter \textit{has an open basis} if
for every $w$ there is an open $v$ in the filter with $v
\subseteq w$. As in general there need not be a point
contained in the intersections of the sets in a filter,
we are interested in a condition on filters that requires
their sets to contain a common point. Our condition that a
filter meets a semi-topogenous order will imply that this
filter has also an open basis, while a compelteness
condition will ensure that each filter that meets the order
has a nonempty intersection.

\begin{definition} 
Let $X$ a space with a semi-topogenous order $\sqsubset$ generating
its topology. A filter $U$ on $X$ \textit{meets}  $\sqsubset$ if for
every $w \in U$ there is a $v \in U$ with $v \sqsubset w$. $X$ is
\textit{complete for $\sqsubset$} if for every filter $U$ in the space
$X$ which meets $\sqsubset$ there is a point $x$ with $x \in u$ for all
$u \in U$.
\end{definition}

%Let $X$ be a $T_1$ topological space. The next result shows
%that the existence of a semi-topogenous order that generates
%the topology for $X$, such that $X$ is complete for the
%order, is sufficient for $X$ to be homeomorphic to an MF space.

\begin{theorem} \label{s7t1} 
Let $X$ be a $T_1$ space with a semi-topogenous order $\sqsubset$ 
generating its topology such that $X$ is complete for
$\sqsubset$. Then $X$ is homeomorphic to an MF space. 
\end{theorem}

\begin{proof}
  Let $P$ consist of the nonempty open subsets of $X$ and 
  let $p \prec q$ if and only if $p \neq q$ and 
  $p \sqsubset q$. The relation $\prec$ is obviously transitive and
  antireflexive, and so it makes $P$ into a poset.

  For each $x \in X$ let $U_x$ be the set of all $p \in P$
  with $\{x\} \sqsubset p$. If $p,q \in U_x$ then the open
  kernel $u$ of $p \cap q$ contains $x$ and thus there is an open
  $r \sqsubset u$ with $x \in r$. As the open kernel of $r$
  again contains $x$, $\{x\} \sqsubset r$. So $r \in U_x$,
  $r \sqsubset p$ and $r \sqsubset q$. Thus $U_x$ is a
  filter on $P$.
  
%  Assume now that $V$ is a maximal filter in $\MF(P)$; then
%  $V$ is also a filter of sets which meets $\sqsubset$: If
  If $V$ is a maximal filter on $\MF(P)$ then $V$ also meets $\sqsubset$.
  If $v$ generates $V$ then $v$ is open (by definition) and not empty. For
  every $x \in v$ there is a $w \sqsubset v$ with $x \in w$;
  by maximality $w = v$. Thus $v \sqsubset v$ and every $w
  \subseteq X$ with $v \subseteq w$ satisfies $v \sqsubset
  w$ and $w \in V$. If there is no single element generating
  $V$ then there is, for every $v, w \in V$, some $u \in V$ with
  $u \prec v$ and $u \prec w$. Then it follows that $u
  \sqsubset v$ and $u \sqsubset w$. Furthermore, there is an
  $t \prec u$ with $t \in V$; then it follows that $t
  \sqsubset v \cap w$. So $V$ contains all supersets of $v
  \cap w$ and thus $V$ is a filter. Furthermore, $V$ meets
  $\sqsubset$.

  This means, by assumption, that there is a point $x$ contained in all
  sets of $V$. Thus $V \subseteq U_x$ and by the maximality
  of $V$, $V = U_x$.  Therefore, every filter $U_y$ is
  contained in a filter $U_x$ which is maximal. Due to the
  $T_1$ property, $y = x$; otherwise $U_y$ would contain a
  $p$ with $x \notin p$ in contradiction to the fact that
  $U_y \subseteq U_x$.
 
  This shows that the mapping $\phi\colon x \mapsto U_x$ is
  a bijection from $X$ to the maximal filters on $P$. To see
  the $\phi$ is open and continuous, first note that if $y
  \in X$ and $U$ is an open set, then $y \in U$ if and only
  if $\{y\} \sqsubset U$. To see this, fix $y \in X$ and any
  open $U$ such that $y \in U$, which means then $\{y\}
  \subseteq U$.  Then, because $\sqsubset$ generates the
  topology and $y$ is trivially in the open kernel of $U$,
  there is some $W \sqsubseteq U$ with $y \in W$. This means
  $\{y\} \subseteq W \sqsubset U$, which means $\{y\}
  \sqsubset U$ by the definition of semi-topogenous orders. The
  converse direction of the equivalence follows directly from the
  definition of a semi-topogenous order.

  Now, to see that $\phi$ is open and continuous, note that for any
  point $x \in X$ and any open set $U$, we have
  \[
  x \in U \Leftrightarrow \{x\} \sqsubset U \Leftrightarrow U \in
  \phi(x) \Leftrightarrow \phi(x) \in N_p,
  \]
  where $N_p$ is the basic open subset of $\MF(P)$ corresponding to~$U$.
\end{proof}

\noindent
We do not know whether every MF space has a semi-topogenous
order satisfying the hypotheses of the previous theorem.
We have established the following partial result.

\begin{theorem}
 If $X = \MF(P)$ and $P$ satisfies
  \begin{equation}\label{eqntocond}
    \forall p,q,r \, [p \prec q \land N_q \subseteq N_r
    \Rightarrow p \prec r]
  \end{equation}
  then there is a semi-topogenous order $\sqsubset$
  generating the topology of $X$ such that $X$ is complete
  for $\sqsubset$.
\end{theorem}

\begin{proof}
  For any $v,w \subset X$, let $v \sqsubset w$ if either $v=\varnothing$,
  $w=X$,  there is is an open atom $u$ with $v \subseteq u \subseteq w$
  or there are $p,q \in P$ with $v \subseteq N_p$, $p \prec q$ and
  $N_q \subseteq w$. Note that $N_p \subseteq N_q$ in the last case.

  It follows directly from definitions and the present
  assumptions that $\sqsubset$ is a semi-topogenous order.
  We must show that that $\sqsubset$ generates the topology
  of $X$.  Let $w$ be an open set and $x$ be a point in $w$.
  There is an open set $N_q$ with $\{x\} \subseteq N_q
  \subseteq w$. In the case that $\{x\} = N_q$, $N_q
  \sqsubset w$. In the case that $\{x\} \neq N_q$ there is a
  further $p \prec q$ with $x \in N_p$: The reason is that
  given an $y \in N_q\setminus\{x\}$, the maximal filter
  $U_x$ belonging to $x$ must contain a $p \prec q$ which
  does not contain $y$ by the $T_1$ axiom. Then $\{x\}
  \subseteq N_p \sqsubset w$. So $w$ is the union of all $v$
  with $v \sqsubset w$.

  Now let $W$ be a filter in the topological space $X$ which
  meets $\sqsubset$. If $W$ contains an $r$ such that $N_r$ is atomic,
  that is, a singleton $\{x\}$, then every $u \in W$ contains $x$
  since otherwise $N_r \cap u = \varnothing$ in contradiction
  to $W$ being a filter.

  If $W$ does not contain an $r$ such that $N_r$ is atomic,
  then let $V$ be the set of all $p \in P$ such that $N_p
  \in W$. Given any $p,q \in V$ there is an $u$ such that $u
  \sqsubset N_p \cap N_q$. So there are $r,t$ with $u
  \subseteq N_r$, $N_t \subseteq N_p \cap N_q$ and $r \prec
  t$. It follows that $r \prec p$ and $r \prec q$. Thus $V$
  is the basis of a filter on $P$; this filter is contained
  in a maximal filter on $P$ which is of the form $U_x$ for
  some point $x \in X$. This $x$ is then in $N_p$ for all
  $p \in V$. Let $u \in W$. As $W$ meets $\sqsubset$, there
  is a $p \in V$ with $N_p \subseteq u$. It follows that $x
  \in N_p$ and $x \in u$. So $x$ is a common point of the
  sets in $W$.
\end{proof}

\noindent
The posets constructed in Theorems~\ref{s2t1}
and~\ref{s2t2} satisfy condition~(\ref{eqntocond}) and thus
they are examples of a poset space that is complete for a
semi-topogenous order generating its topology.

\begin{example}
  For every complete metric space and every locally compact Hausdorff space
  there exists a semi-topogenous order $\sqsubset$ which generates the
  topology of $X$ and for which $X$ is complete.
\end{example}

\begin{remark}
  Assume that $X$ is a space which is complete for a
  semi-topogenous order generating its topology. Then one
  can not only show that $X$ is homeomorphic to an MF space
  but also that the winning strategy for player~II is quite
  easy to obtain. Given any open set $u$ and any point $x \in u$ by
  player~\I, player~\II only has to choose an open $v$ with
  $\{x\} \subseteq v \sqsubset u$. It does not matter which
  $v$ with this condition is chosen and the history of the
  game can be ignored. The result of the construction will
  be, at the end of the game, a basis for a filter which
  meets $\sqsubset$ and thus this filter has a common point.

  This shows that the ``neighbourhood spaces'' that we consider
  here satisfy a restricted version of the strong Choquet
  property. The intuition behind this restriction is that
  one wishes to study non-second-countable spaces by
  considering ``transfinite games.'' The role of player~\I
  is replaced by considering filters instead of descending
  sequences, and the winning strategy of player~\II is
  reduced to a neighbourhood relation $\sqsubset$ which
  could be interpreted as saying that if $\{x\} \subseteq v
  \sqsubset u$ then $v$ is a good move for player~\II.

  Indeed, the notion of completeness of spaces with respect
  to a semi-topogenous order $\sqsubset$ is based on this
  idea: Let the strategy of player~\II be just to follow
  $\sqsubset$ and let player~\I build a filter $U$ such that
  for every $w \in U$ there is a $v \in U$ which player~\II
  might have chosen as a response to $w$, that is, a $v
  \sqsubset w$; then the intersection of all $u \in U$ is
  not empty.
\end{remark}

\section{Cardinality of poset spaces}\label{sec6}

\noindent
In this section, we establish perfect set theorems for
countably based Hausdorff poset spaces. These theorems show
that these spaces are either countable or have the
cardinality $2^{\aleph_0}$ of the continuum.

\begin{theorem}\label{s8t1}
  Any countably based Hausdorff poset space has either
  countably many points or has cardinality $2^{\aleph_0}$.
\end{theorem}

\begin{corollary}
  Any countably based Hausdorff poset space has either
  countably many points or contains a perfect closed set.
\end{corollary}

\begin{proof}
  Any second-countable Hausdorff space of cardinality
  $2^{\aleph_0}$ contains a perfect closed set. The complement
  of the perfect closed set is the union of all the basic
  open sets from a fixed countable basis which contain fewer
  than $2^{\aleph_0}$ points. 
\end{proof}

\noindent
To prove Theorem~\ref{s8t1}, we introduce a class of
Gale--Stewart games. These games are inspired by the
$*$-games in descriptive set theory
(as described in~\cite{Kechris-CDST}). For each poset $P$ we define a
game which we call the \define{poset star game on $P$}.
There are two players. The play proceeds in stages numbered
$0,1,2,\dots$. At stage $t$, player~\I plays a pair
$\langle p^t_1,p^t_2\rangle \in P \times P$. Player~\II
plays a number $n_t \in \{1,2\}$. Player~\I wins the game
if the following conditions hold for all $t$.
\begin{enumerate}
\item $p^t_1 \perp p^t_2$.
\item $p^{t+1}_1 \preceq p^t_{n_t}$ and  $p^{t+1}_2 \preceq p^t_{n_t}$.
\end{enumerate}
Player~\II wins if player~\I does not win; there are no
ties.

A \define{strategy} for a player is a function that tells
the player what to do at any possible move of the game. The
strategy is a \define{winning strategy} if the player will
win any play of the game in which the player uses the
strategy to choose every move. It is impossible for both
players to have a winning strategy for the same game.

\begin{lemma}\label{lem:stardet}
  Let $P$ be a poset. Either player~\I or player~\II has a
  winning strategy for the poset star game on $P$.
\end{lemma}

\begin{proof}
  The set of infinite plays of the poset star game on $P$
  that are winning for player~\I is closed in the space of all
  possible plays of the game. (This space is the space of
  infinite sequences of moves; the set of moves is given the
  discrete topology and the space of infinite plays carries
  the product topology). The proof follows from a
  theorem of Gale and Stewart known as closed determinacy. 
\end{proof}

\begin{lemma}\label{lem:stardeti}
  Suppose that $X$ is a Hausdorff poset space based on a countable
  poset $P$ and player~\I has a winning strategy
  for the poset star game on~$P$. Then $X$ has cardinality
  $2^{\aleph_0}$.
\end{lemma}

\begin{proof}
  It suffices to prove the result for $\MF(P)$, which is a
  subset of $\UF(P)$. Let $s_\I$ be a winning strategy for
  player~\I and let $f \in \{1,2\}^\setN$. Consider the play in
  which player~\I follows $s_\I$ while player~\II uses $f$
  as a guide; that is, player~\II plays $f(n)$ at stage
  $2n$. Because $s_\I$ is a winning strategy for player~\I,
  this play determines a descending sequence $F(f)$ of
  elements of $P$. This sequence extends to a maximal
  filter. For distinct $f,g \in \{0,1\}^\setN$ the sequences
  $F(f)$ and $F(g)$ contain incompatible elements and thus
  cannot extend to the same filter. Therefore the space
  $\MF(P)$ has cardinality $2^{\aleph_0}$.
\end{proof}

\begin{lemma}\label{lem:stardetii}
  Let $X$ be a countably based Hausdorff poset space based
  on the poset $P$. If player~\II has a winning strategy
  for the poset star game on $P$ then $X$ is countable.
\end{lemma}

\begin{proof}
  Let $s_{\II}$ be a winning strategy for player~\II. We say that
  a finite play $\sigma$ of length $2k$ is \textit{compatible
    with $s_{\II}$} if $s_{\II}(\sigma[2i + 1]) = \sigma(2i +2)$
  whenever $2i + 2 \leq k$. We say that a play $\sigma$ of even
  length is a \textit{good play} for a point $x$ if $\sigma$ is
  compatible with $s_{\II}$ and $x$ is in the open set
  chosen by player~\II in the last move of $\sigma$. A good play for
  $x$ is a \textit{maximal play} if it cannot be extended to
  a longer good play for $x$; this means that no matter what
  pair of disjoint open sets player~\I plays, $s_{\II}$ will direct
  player~\II to choose an open set not containing $x$.
  
  If player~\II has a winning strategy then every point $x$ has a
  maximal play. Note that the empty play is trivially a
  good play for $x$. If every good play for $x$ could be
  extended to a larger good play for $x$, then it would be
  possible for player~\I to win the game by always leaving
  the game in a position that is good for $x$. This play of
  the game would follow $s_{\II}$, a winning strategy for
  player~\II, which is a contradiction.

  If $\sigma$ is a good play for two points $x$ and $y$ then
  $\sigma$ is not a maximal play for both $x$ and $y$. For
  player~\I could play $\langle U_1, U_2\rangle$ in response to
  $\sigma$, where $x \in U_1$, $y \in U_2$ and $U_1 \cap
  U_2 = \varnothing$. Here we are using the assumption that
  the topology of $X$ is Hausdorff.

  We have now shown that every point in the $X$ has a
  maximal play and that no play is maximal for two points.
  Since the set of maximal plays is countable, this implies
  that the set of points in $X$ is countable.
\end{proof}

\noindent
We remark that the statement ``Every closed subset of a
countably based Hausdorff MF space is either countable or
has a perfect closed subset'' is independent of ZFC set
theory; this result is established in~\cite{Mummert-MFspaces}.

\section*{Acknowledgments}

\noindent
We would like to thank the Institute for Mathematical
Sciences at the National University of Singapore for
organizing the wonderful Computational Prospects of Infinity
workshop in 2005 which made this work possible. We also
thank Steffen Lempp and Sasha Rubin for thoughtful comments.
We would like to thank Jimmie Lawson and Ralph Koppermann
for their helpful comments on the domain theory results in
Section~\ref{sec5a}.

Some of the results presented here appeared in the first
author's PhD thesis~\cite{Mummert-thesis}, supervised by
Stephen Simpson at the Pennsylvania State University.

%\newpage

\noindent\textbf{Carl Mummert}\\
Department of Mathematics\\
Marshall University \\
One John Marshall Drive\\
Huntington, WV 25755\\
USA\\
\mbox{}

\noindent\textbf{Frank Stephan}\\
School of Computing and Department of Mathematics\\
National University of Singapore\\
Singapore 117543\\
Republic of Singapore

\end{document}